\documentclass[12pt]{amsart}
\usepackage{amsmath}
\begin{document}
\title{The Erd\"{o}s and Campbell-Staton Conjectures about Square Packing}
\author{iwan Praton}
\address{Department of Mathematics, Franklin and Marshall College, Lancaster PA 17604}
\email{iwan.praton@fandm.edu}
\date{}

\begin{abstract}
Put $n$ open non-overlapping squares inside a unit square, and let $f(n)$ denote
the maximum possible value of the sum of the side lengths of the $n$ squares. 
Campbell and Staton, building on a question of Erd\"{o}s, conjectured that 
$f(k^2+2c+1)=k+c/k$, where $c$ is any integer and $k\geq |c|$. We show that
if this conjecture is true for one value of $c$, then it is true for all values of $c$.
\end{abstract}

\maketitle
Let $n$ be a positive integer. Put $n$ open non-overlapping squares inside a unit 
square. If $e_1,e_2,\dotsc, e_n$ are the side lengths of the squares, then define
$f(n) = \max \{e_1+e_2+\cdots +e_n\}$. It is not hard to show that $f(k^2)=k$; for a
proof see~\cite{1}. In the beginning of~\cite{2}, it is mentioned that $f(k^2+1)=k$
is an old conjecture of Erd\"{o}s. Apparently not much work was done on this
problem---even~\cite{2} is mostly about packing unit squares in a larger square---but
in~\cite{1}, Campbell and Staton revived the question, provided some nice lower bounds
for $f(n)$, and conjectured that their lower bounds are actually exact. 

In this short note we show that the Erd\"{o}s conjecture implies the Campbell-Staton
conjectures. Of course, we should describe the conjectures first. Let $c$ be
any integer. Then the Erd\"{o}s and Campbell-Staton conjectures can be described
succinctly as follows:
\[
f(k^2+2c+1)=k+c/k \text{ for all $k\geq |c|$}.\leqno{(*)}
\]
When $c=0$ we get the Erd\"{o}s conjecture; Campbell and Staton put forth
all the other values of $c$. They also showed, by explicit construction, 
that $f(k^2+2c+1)\geq k+c/k$ for all $k\geq |c|$, so all that remains to 
prove the conjecture is to provide the appropriate upper bound. 
Note that if $n$ is any nonsquare integer, then $n$
is an odd integer away from a neighboring square, so $(*)$ provides a complete
description of the function $f$.

We'll show that if $(*)$ is true for one value of $c$, then it is true for all
values of $c$. In particular, as stated above, the truth of the Erd\"{o}s conjecture
implies the truth of the Campbell-Staton conjectures. Let $P(c)$ denote 
the statement $(*)$. It clearly suffices to show that
\[
P(c-1)\implies P(c) \text{ and } P(c+1)\implies P(c);
\]
we'll show this after deriving a general upper bound for $f(n)$.

Take a unit square and divide it into the standard $b\times b$
grid of squares, each with side length $1/b$. 
Remove an $a\times a$ subsquare, and replace it with an optimal
configuration
of $n$ squares, shrunk by a factor of $b/a$
so that it fits inside the $a\times a$ space.
We now have a configuration of
$b^2-a^2+n$ squares inside the unit square. The sum of the side lengths
of these squares is $af(n)/b+(b^2-a^2)/b$. This is at most $f(b^2-a^2+n)$,
so we have $af(n)/b+(b^2-a^2)/b\leq f(b^2-a^2+n)$. Solving for $f(n)$
gives us an upper bound for $f(n)$:
\[
f(n)\leq a-\frac{b^2}{a} + \frac{b}{a}f(b^2-a^2+n).\leqno{(**)}
\]

We now prove that $P(c-1)\implies P(c)$. Let $k\geq |c|$, and
put $n=k^2+2c+1$, $a=k-1$, $b=k$ in $(**)$. 
Then $b^2-a^2+A = 
2k-1+k^2+2c+1=(k+1)^2+(2c-1)=(k+1)^2+2(c-1)+1$. Note that 
 $k+1\geq |c-1|$, so we can apply $P(c-1)$, i.e., 
$f(b^2-a^2+n)=k+1+(c-1)/(k+1)=k+(k+c)/(k+1)$. Thus $(**)$ becomes
\begin{align*}
f(k^2+2c+1)&\leq k-1-\frac{k^2}{k-1} + \frac{k}{k-1}\left(k+\frac{k+c}{k+1}\right)\\
&=k+\frac{c}{k} + \frac{k+c}{k(k^2-1)}.
\end{align*}
after some straightforward calculations.

This is not quite what we want---it's too big by $(k+c)/(k(k^2-1))$. 
But we can use $(**)$ again. Let $n=k^2+(2c+1)$
as before, but now let $a=k$ and keep $b$ arbitrary. Then $b^2-a^2+n=
b^2+(2c+1)$. We get
\begin{align*}
f(k^2+2c+1) &\leq k-\frac{b^2}{k}+\frac{b}{k}f(b^2+(2c+1))\\
&\leq k-\frac{b^2}{k}+\frac{b}{k}\left(b+\frac{c}{b}+\frac{b+c}{b(b^2-1)}\right)\\
&=k+\frac{c}{k}+\frac{b+c}{k(b^2-1)}.
\end{align*}
Now let $b\rightarrow \infty$. We get
\[
f(k^2+2c+1)\leq k+c/k,
\]
which, in conjunction with the Campbell-Staton lower bounds, gives us what we want.

The implication $P(c+1)\implies P(c)$ is proved similarly, but we use $n=k^2+2c+1$,
$a=k+1$, and $b=k+2$ in the first step. Details are left to the reader.

As a final remark, note that in order to  prove the conjectures completely,
it is enough to show that $f(k^2+2c+1)=k+c/k+\epsilon(k)$,
where $k\epsilon(k)\rightarrow 0$ as $k\rightarrow \infty$. To do this it is probably
necessary to consider the detailed placements of the small squares in
the unit square.

\end{document}